\documentclass[11pt]{article}

\usepackage[english]{babel}
\usepackage{graphicx}
\usepackage{subcaption}

\newcommand{\R}{\mathbb{R}}

\newcommand{\Z}{\mathbb{Z}}
\newcommand{\C}{\mathbb{C}}

\usepackage{siunitx}

\usepackage[margin=2.5cm]{geometry}
 \usepackage{amsmath, amssymb}
\usepackage{amsthm}

\usepackage{xcolor}
\usepackage[colorlinks = true,
            linkcolor = blue,
            urlcolor  = blue,
            citecolor = blue,
            anchorcolor = blue]{hyperref}
      
\newtheorem{theorem}{Theorem}[section]
\newtheorem{remark}[theorem]{Remark}

\usepackage{titling}
\newcommand*\samethanks[1][\value{footnote}]{\footnotemark[#1]}

\begin{document}

\title{An analytical solution for vertical infiltration in bounded profiles}

\author{I. Argyrokastritis\thanks
  {Department of Natural Resources Development and Agricultural Engineering, Agricultural University of Athens, Greece}, \, K. Kalimeris\thanks{Mathematics Research Center, Academy of Athens, Greece} \, and L. Mindrinos\samethanks[1]}

\thanksmarkseries{arabic}

\maketitle

\begin{abstract}
In this study, we derive an analytical solution to address the problem of one-dimensional vertical infiltration within bounded profiles. We consider the Richards equation together with various boundary conditions, simulating different scenarios of water application onto the surface of a homogeneous and bounded medium. To solve the corresponding initial boundary value problem over a finite interval, we apply the unified transform, commonly known as the Fokas method. Through this methodology, we obtain an integral representation that can be efficiently and directly computed numerically, yielding a convergent scheme.
\end{abstract}

\section{Introduction}\label{intro}

The term infiltration defines the process of water flow into the soil through its surface. Infiltration may be one-, two-, or three-dimensional. Especially vertical infiltration of water into the soil mass is a process of great importance in soil hydrology because its rate regulates the amount of the water that penetrates the soil surface, to increase the water content of the soil profile, and the amount of runoff over the soil surface during heavy rains. There are many soil factors affecting the infiltration process, like the characteristics of the soil, the condition of the soil surface and the deficit of soil moisture. However the main factor is the hydraulic conductivity and its dependence on the water content. 

If water is applied by flooding to the surface of a homogeneous and deep soil with uniform initial moisture, saturated conditions are created at the surface in time zero and water enters the soil through its surface at a rate called infiltration rate. The rate of water infiltration decreases over time. This reduction is due to the reduction of the hydraulic gradient on the surface but also to other factors such as sealing of the soil surface. After enough time has passed, the infiltration rate stabilizes at a final value which is equal to the saturated hydraulic conductivity. 

If water is applied to the soil surface with a constant rate (e.g. rainfall), the infiltration rate in the initial stages is equal to the rate of supply and depends on the provision of application and not on the conditions and properties of the soil. In this case until some time called time to ponding, all water applied is infiltrated until the soil surface is saturated. After that time part of the water applied does not penetrate the flooded soil surface and flows superficially as surface runoff. This means that the actual infiltration rate is smaller than rainfall intensity and later it takes its final value, which is equal to the saturated hydraulic conductivity.

The mathematical treatment of the phenomenon of water infiltration into the soil mass includes the definition of initial and boundary conditions that characterize the soil. The hydraulic characteristic functions of the soil must also be known as they contribute to the solution of Richards equation, which is the main equation that describes the infiltration process under the initial and boundary conditions imposed. 

Mathematical procedures for the solution of infiltration problems, are either numerical or analytical.  Numerical solutions are obtained by using various finite difference or finite element numerical schemes and their modifications and/or improvements.  For an overview, we refer to the initial works \cite{Haverkamp77, Rubin63} and the review papers \cite{Farthing17, feddes1988, ross1990, van00} and the references therein. On the other hand, analytical solutions do not exist always. We can derive such only for, in a way, trivial boundary conditions of unsteady flow in soils, which is defined by simple hydraulic functions and for steady flow processes in homogeneous and layered soil profiles. 

Analytical solutions developed for elementary soil hydrological processes, like infiltration, prove suitable for soils characterized by straightforward hydraulic functions. Noteworthy examples include Carslaw and Jaeger's analytical solutions for infiltration with flooding \cite{carslaw1959} and Braester's solutions for rainfall \cite{Bra73}, particularly applicable when diffusivity and conductivity remain constant. In cases where diffusivity is constant and conductivity varies as a function of water content, Philip provides analytical solutions for linear and quadratic dependencies \cite{philip66, philip1974}. It is crucial to note, however, that all the aforementioned solutions are regarded as ``approximate" analytical expressions due to the appearance of the complementary error function in their formulas. 

Such solutions, even if they are not applicable directly to field conditions have a great advantage. They lead to full understanding of the physics of the infiltration phenomenon, and they provide estimates of deviations caused by any reason e.g. deviations caused by possible alteration of a boundary condition. Another advantage of these solutions is that they allow the quantification of errors of approximate and numerical procedures used for the description of infiltration. 

In the current work we present the mathematical modelling of the vertical infiltration in bounded profiles, by formulating an initial boundary value problem in the finite interval $(0,L)$. For example, the physical condition of prohibiting  any infiltration at the bottom of that interval can be modelled by the boundary condition which fixes the water content at $x=L$; this is to be contrasted with the half-line modelling, equivalently $L\to\infty$.  Then, we derived the solution of that problem in the form of an integral representation involving transforms of the initial and boundary conditions.

The derivation of this solution representation is based on the Fokas method for solving linear initial-boundary value problems (IBVPs). This method, which is also known as the unified transform, was introduced in 1997 by Fokas \cite{F97}, for solving IBVPs for nonlinear integrable partial differential equations (PDEs). Later, it was realised that it produced effective analytical and numerical solutions for linear PDEs \cite{F02}; for an overview of the method we refer to \cite{F08, FK22}. The last two decades hundreds of works have been published, using and extending the Fokas method to substantially different directions, including  control theory \cite{KO20, KOD23, OK23} and regularity results \cite{Him20, OY19, BFO20}; we consider these, and references therein important directions for the development of the current results to other interesting directions.  In general, the Fokas method yields solutions of the IBVPs as integral representations, with exponentially decaying integrands, involving the Fourier (spectral) transforms of the initial and boundary conditions.  This is important for analytical purposes \cite{AMO24, F17}, as well as for numerical purposes \cite{KOD23, barros19, CFH19}. In the current work, we derive the solution for general initial and boundary conditions and we illustrate the distribution of the water content for various physical set-ups exploiting the numerical advantages provided by the Fokas method.

The paper is organized as follows: In  \autoref{sec_formulation} we formulate mathematically the problem of vertical infiltration in a bounded domain as an initial boundary value problem. Then, we obtain its analytical solution using the Fokas method. Its integral representation in the complex plane contains transforms of the initial and the known boundary functions. In \autoref{sec_examples} we consider three different examples modelling  flooding and rainfall by varying the initial and the boundary functions. The numerical implementation and the comparison of our solution with well-known and mostly used formulas is performed in \autoref{sec_numerics}.

\section{Problem formulation}\label{sec_formulation}

The vertical infiltration in an one-dimensional homogeneous porous medium (for example soil) is modelled through the Richards equation
\begin{equation}\label{eq_rich}
\frac{\partial \theta}{\partial t} = \frac{\partial}{\partial x} \left(K(\theta) \frac{\partial h}{\partial x} \right) - \frac{\partial K (\theta)}{\partial x}, \quad \mbox{for  }x>0, \,\, t>0,
\end{equation}
where $\theta$ is the water content, $K$ is the hydraulic conductivity and $h$ is the water pressure.  Using the relation (first introduced in \cite{childs})
\begin{equation}
D(\theta) = K (\theta) \frac{d h}{d\theta},
\end{equation}
where $D$ describes the water diffusivity, we obtain
\begin{equation}\label{diffusion}
\frac{\partial \theta}{\partial t} = \frac{\partial}{\partial x} \left(D(\theta) \frac{\partial \theta}{\partial x} \right) - \frac{\partial K (\theta)}{\partial x},
\end{equation}
a non-linear diffusion-type PDE.

The non-linear dependence of $D$ and $K$ on $\theta$ prevent us from deriving an analytical  solution. Thus, a commonly done simplification is to consider specific models where the relations between the conductivity and diffusion coefficients and $\theta$ are well described.  A typical assumption is to consider power law relations of the form 
\begin{equation}\label{power_law}
D(\theta) = D_0 \theta^n, \quad K(\theta) = K_0 \theta^k, \quad \mbox{for  } n,k=1,2,\ldots 
\end{equation}
for constant coefficients $D_0$ and $K_0.$ The case $n=0$ and $k=1,$ was first examined in \cite{philip1957} and an analytical solution was derived. Later, Philip considered also the case $n=0$ and $k=2$ in \cite{philip1974}. More recently, an exact solution was derived for the general case of $n=k-1,$ for $k=2,3,\ldots$ \cite{hayek16}.  In a different context, De Barros et al. \cite{barros19} considered the case $n=0,$ and $k=1,$ to derive an analytical solution on the half-line using the unified transform. 

In this work, we are interested in solving \eqref{diffusion}, considering \eqref{power_law} for $n=0$ and $k=1,$  in a bounded domain using the unified transform in order to derive an analytical solution that leads to convergent numerical schemes. This can be seen as a first step to model the water flow problem in layered soils. 

To be more precise, we consider the following initial boundary value problem:
\begin{subequations}\label{bvp}
\begin{alignat}{3}
\frac{\partial \theta}{\partial t} + K_0 \frac{\partial \theta}{\partial x} &= D_0 \frac{\partial^2 \theta}{\partial x^2},  \quad && 0<x <L, \, t>0,  \label{bvp1}\\
\theta (x,0) &= \theta_0 (x), \quad && 0 <x <L, \label{bvp2}\\ 
\theta (0,t) = f(t), \,\, \theta (L,t) &= g(t), \quad &&t>0, \label{bvp3}
\end{alignat}
\end{subequations}
for general spatial-dependent initial function $\theta_0$ and time-dependent boundary functions $f$ and $g.$ Specific functions will be considered later for the numerical implementation, modelling either rainfall infiltration (varying $f$) or infiltration with constant head - flooding at the soil surface (constant $f$).

\subsection{Analytical solution}

In the section, we derive an analytical solution of \eqref{bvp} using the Fokas method. We multiple \eqref{bvp1} with $e^{-i \lambda x}$ and we integrate with respect to $x,$ to obtain
\begin{equation}\label{integral}
\frac{\partial}{\partial t} \int_0^L e^{-i \lambda x} \theta (x,t) dx + K_0 \int_0^L e^{-i \lambda x} \frac{\partial \theta}{\partial x} (x,t) dx  = D_0 \int_0^L e^{-i \lambda x} \frac{\partial^2 \theta}{\partial x^2} (x,t) dx.
\end{equation}
The second term using integration by parts together with \eqref{bvp3} reads
\begin{equation}\label{parts1}
\int_0^L e^{-i \lambda x} \frac{\partial \theta}{\partial x} (x,t) dx = e^{-i \lambda L} g(t) - f(t) + i \lambda \hat\theta (\lambda,t),
\end{equation}
where we have defined
\begin{equation}\label{fourier}
\hat\theta (\lambda,t) =  \int_0^L e^{-i \lambda x} \theta (x,t) dx.
\end{equation}
Similarly, the last term in \eqref{integral} takes the form
\begin{equation}\label{parts2}
\int_0^L e^{-i \lambda x} \frac{\partial^2 \theta}{\partial x^2} (x,t) dx = e^{-i \lambda L} \frac{\partial \theta}{\partial x} (L,t)-\frac{\partial \theta}{\partial x} (0,t)  
+ i \lambda \left(e^{-i \lambda L} g(t) - f(t) \right)-\lambda^2 \hat\theta (\lambda,t).
\end{equation}
We substitute \eqref{parts1} and \eqref{parts2} in \eqref{integral} to get
\begin{equation}\label{integral2}
\frac{\partial}{\partial t}\hat\theta (\lambda,t) + \omega (\lambda) \hat\theta (\lambda,t) = e^{-i \lambda L}  (i \lambda D_0 -K_0) g(t) - (i \lambda D_0 -K_0) f(t ) + D_0 \left(e^{-i \lambda L} \frac{\partial \theta}{\partial x} (L,t)-\frac{\partial \theta}{\partial x} (0,t)\right),
\end{equation}
where $\omega (\lambda) =  D_0  \lambda^2 + i K_0 \lambda. $ We multiple both sides of \eqref{integral2} with $e^{\omega(\lambda)t}$ and we integrate with respect to $t,$ resulting in
\begin{equation}\label{integral3}
\begin{aligned}
e^{\omega (\lambda) t}\hat\theta (\lambda,t) -  \hat\theta (\lambda,0) &=
\int_0^t e^{\omega (\lambda) \tau} \left[
  e^{-i \lambda L}  (i \lambda D_0 -K_0) g(\tau) - (i \lambda D_0 -K_0) f(\tau )
 \right. \\
   &\phantom{=} \left. + D_0 \left(e^{-i \lambda L} \frac{\partial \theta}{\partial x} (L,\tau)-\frac{\partial \theta}{\partial x} (0,\tau)\right) \right] d\tau.
\end{aligned}
\end{equation}
We define
\begin{equation}
\begin{aligned}
\tilde f_0 (\lambda,t) &= \int_0^t e^{\lambda \tau} f(\tau )d\tau,\\
\tilde f_1 (\lambda,t) &= \int_0^t e^{\lambda \tau} \theta_x (0,\tau )d\tau, \\
\tilde g_0 (\lambda,t) &= \int_0^t e^{\lambda \tau} g(\tau )d\tau, \\
\tilde g_1 (\lambda,t) &= \int_0^t e^{\lambda \tau} \theta_x (L,\tau )d\tau,
\end{aligned}
\end{equation}
and we rewrite \eqref{integral3} in the compact form
\begin{equation}\label{global}
\begin{aligned}
\hat\theta (\lambda,t)  &= e^{-\omega (\lambda) t} \left[
  \hat \theta_0 (\lambda) + e^{-i \lambda L}  (i \lambda D_0 -K_0) \tilde g_0 (\omega (\lambda),t) - (i \lambda D_0 -K_0) \tilde f_0 ( \omega (\lambda), t) \right. \\
  &\phantom{=} \left. + D_0 \left(e^{-i \lambda L}  \tilde g_1 (\omega (\lambda),t)  - \tilde f_1 (\omega (\lambda),t)\right) \right].
  \end{aligned}
\end{equation}

This is the global relation between the initial and boundary conditions. Following \cite{barros19}, the application of the inverse Fourier transform results in the integral representation of the solution
\begin{equation}\label{repre1}
\begin{aligned}
\theta (x,t) &= \frac1{2\pi} \int_{\R} e^{i \lambda x - \omega (\lambda)t} \hat \theta_0 (\lambda) d\lambda - \frac1{2\pi} \int_{\partial D_+} e^{i \lambda x - \omega (\lambda)t} \left[(i \lambda D_0 -K_0) \tilde f_0 ( \omega (\lambda), t) + D_0 \tilde f_1 (\omega (\lambda),t) \right]d\lambda \\
&\phantom{=}- \frac1{2\pi} \int_{\partial D_-} e^{-i \lambda (L-x) - \omega (\lambda)t} \left[(i \lambda D_0 -K_0) \tilde g_0 (\omega (\lambda), t) + D_0 \tilde g_1 (\omega (\lambda),t) \right]d\lambda,
\end{aligned}
\end{equation}
where $\partial D_{\pm} = \big\{ \lambda \in \C: \operatorname{Re}(\omega (\lambda)) = 0, \ \operatorname{Im}(\lambda) \gtrless 0 \big\}$ is the boundary of the domain $D_{\pm} = \big\{ \lambda \in \C: \operatorname{Re}(\omega (\lambda)) < 0, \ \operatorname{Im}(\lambda) \gtrless 0 \big\}.$

In the above expression, two terms ($\tilde f_1$ and $\tilde g_1$) are unknown. In the following, we will eliminate them. We apply the invariant transform $\lambda \mapsto \nu (\lambda) = -\lambda - i\tfrac{K_0}{D_0}$ in the global relation \eqref{global}. 
This transform satisfies $\omega (\lambda) =\omega (\nu(\lambda)).$ We get
 \begin{equation}\label{global2}
\begin{aligned}
\hat\theta (\nu(\lambda),t)  &= e^{-\omega (\lambda) t} \left[
  \hat \theta_0 (\nu(\lambda)) + e^{-i \nu(\lambda) L}  (i \nu(\lambda) D_0 -K_0) \tilde g_0 (\omega (\lambda),t) - (i \nu(\lambda) D_0 -K_0) \tilde f_0 ( \omega (\lambda), t) \right. \\
  &\phantom{=} \left. + D_0 \left(e^{-i \nu(\lambda) L}  \tilde g_1 (\omega (\lambda),t)  - \tilde f_1 (\omega (\lambda),t)\right) \right].
  \end{aligned}
\end{equation}
The functions $\tilde f_1$ and $\tilde g_1$ are the solutions of the system of equations \eqref{global} and \eqref{global2}. We define $\Delta (\lambda) = e^{-i\lambda L}- e^{-i \nu (\lambda) L}$ and 
\begin{equation}\label{funG}
G (\lambda,t) = \hat \theta_0 (\lambda) + e^{-i \lambda L}  (i \lambda D_0 -K_0) \tilde g_0 (\omega (\lambda),t) - (i \lambda D_0 -K_0) \tilde f_0 ( \omega (\lambda), t),
\end{equation}
which is a known function. Then, we get
\begin{equation}\label{sol_system}
\begin{aligned}
\tilde f_1 (\omega (\lambda),t) &= \frac{1}{D_0} \left[ \frac{e^{-i \nu (\lambda) L}}{\Delta (\lambda)} \left(e^{\omega(\lambda)t }\hat\theta (\lambda,t) - G(\lambda,t) \right) -\frac{e^{-i  \lambda L}}{\Delta (\lambda)}\left(e^{\omega(\lambda)t } \hat\theta (\nu(\lambda),t) - G(\nu(\lambda),t) \right)  \right], \\
\tilde g_1 (\omega (\lambda),t) &= \frac{1}{D_0 \Delta (\lambda)} \left[ e^{\omega(\lambda)t }\left(\hat\theta (\lambda,t)-\hat\theta (\nu(\lambda),t)\right) - G(\lambda,t) +G(\nu(\lambda),t) \right].
\end{aligned}
\end{equation}
We substitute \eqref{sol_system} in \eqref{repre1} to derive
\begin{equation}\label{repre2}
\begin{aligned}
\theta (x,t) &= \frac1{2\pi} \int_{\R} e^{i \lambda x - \omega (\lambda)t} \hat \theta_0 (\lambda) d\lambda - \frac1{2\pi} \int_{\partial D_+} e^{i \lambda x - \omega (\lambda)t} (i \lambda D_0 -K_0) \tilde f_0 (\omega (\lambda), t) d\lambda  \\
&\phantom{=}- \frac1{2\pi} \int_{\partial D_+} e^{i \lambda x - \omega (\lambda)t} \frac1{\Delta(\lambda)} \left( e^{-i \lambda L} G(\nu(\lambda),t )- e^{-i \nu (\lambda) L} G(\lambda,t )\right) d\lambda \\
&\phantom{=}- \frac1{2\pi} \int_{\partial D_-} e^{-i \lambda (L-x) - \omega (\lambda)t} (i \lambda D_0 -K_0) \tilde g_0 ( \omega (\lambda), t) d\lambda\\
&\phantom{=}- \frac1{2\pi} \int_{\partial D_-} e^{-i \lambda (L-x) - \omega (\lambda)t} \frac{1}{\Delta (\lambda)} \left(
G(\nu(\lambda),t ) - G(\lambda,t )
 \right) d\lambda \\
 &\phantom{=}- \frac1{2\pi} \int_{\partial D_+} e^{i \lambda x } \frac1{\Delta (\lambda)} \left( e^{-i \nu (\lambda) L} \hat \theta(\lambda,t )- e^{-i \lambda L} \hat\theta (\nu(\lambda),t )\right)d\lambda \\
 &\phantom{=}- \frac1{2\pi} \int_{\partial D_-} e^{-i \lambda (L-x) } \frac{1}{\Delta (\lambda)} \left(
\hat\theta(\lambda,t ) - \hat\theta(\nu(\lambda),t )
 \right) d\lambda .
\end{aligned}
\end{equation}

The last step is to show that the last two terms can be neglected and to simplify the representation of the remaining terms. Following \cite{FokKal22} we recall that on $\partial D_+$ we have $\operatorname{Im}(\lambda) >0,$ thus the term $e^{i\lambda L}$ is bounded and vanishes exponentially for $|\lambda| \rightarrow \infty.$ Then, the asymptotics of the integrand of the second to last term as $\lambda\to\infty$, yield
\begin{align*}
\frac1{\Delta (\lambda)} \left( e^{-i \nu (\lambda) L} \hat \theta(\lambda,t )- e^{-i \lambda L} \hat\theta (\nu(\lambda),t )\right) &\sim \frac1{e^{-i\lambda L}}\left( e^{-i \nu (\lambda) L} \hat \theta(\lambda,t )- e^{-i \lambda L} \hat\theta (\nu(\lambda),t )\right) \\
&= e^{i 2\lambda L -\tfrac{K_0}{D_0}L} \hat \theta(\lambda,t )-  \hat\theta (\nu(\lambda),t ) \\
&= e^{i \lambda L -\tfrac{K_0}{D_0}L} \int_0^L e^{-i\lambda (\xi - L)} \theta(\xi,t ) d\xi -  \hat\theta (\nu(\lambda),t ).
\end{align*} 
In the last equation both the integral, since $\xi-L<0,$ and the last term are bounded, namely $O\left(\frac{1}{\lambda}\right),\ |\lambda|\to\infty$, and analytic on $\partial D_+$, thus from Jordan's lemma we get that the whole integral (second to last term in \eqref{repre2}) vanishes. Similarly, we can show that the last integral in \eqref{repre2} vanishes on $\partial D_-$ where now $\operatorname{Im}(\lambda) <0.$

The second and the third terms in the right-hand side of \eqref{repre2} are combined as follows:
\begin{equation}\label{int_plus1}
- \frac1{2\pi} \int_{\partial D_+} e^{i \lambda x - \omega (\lambda)t} \frac{\alpha_+ (\lambda,t)}{\Delta (\lambda)}
 d\lambda,
\end{equation}
where
\begin{equation}
\begin{aligned}
\alpha_+ (\lambda,t) &=  (i \lambda D_0 -K_0) \Delta (\lambda) \tilde f_0 (\omega (\lambda), t) +  e^{-i \lambda L} G(\nu(\lambda),t )- e^{-i \nu (\lambda) L} G(\lambda,t ) \\
&= (2 i \lambda D_0 -K_0) e^{-i\lambda L } \tilde f_0 (\omega (\lambda), t) - (2 i \lambda D_0 -K_0) e^{-\tfrac{K_0}{D_0} L } \tilde g_0 (\omega (\lambda), t) \\
&\phantom{=}+ e^{-i\lambda L } \hat \theta_0 (\nu (\lambda)) - e^{-i\nu(\lambda) L } \hat \theta_0 ( \lambda).
\end{aligned}
\end{equation}
Similarly, the fourth and the fifth integrals are rewritten as:
\begin{equation}\label{int_minus1}
- \frac1{2\pi} \int_{\partial D_-} e^{-i \lambda (L-x) - \omega (\lambda)t} \frac{\alpha_- (\lambda,t)}{\Delta (\lambda)}
  d\lambda,
\end{equation}
where 
\begin{equation}
\begin{aligned}
\alpha_- (\lambda,t) &= (i \lambda D_0 -K_0) \Delta (\lambda) \tilde g_0 (\omega (\lambda), t) + G(\nu(\lambda),t )-  G(\lambda,t ) \\
&= (2 i \lambda D_0 -K_0)  \tilde f_0 (\omega (\lambda), t) - (2 i \lambda D_0 -K_0) e^{-i \nu (\lambda) L } \tilde g_0 (\omega (\lambda), t) +  \hat \theta_0 (\nu (\lambda)) -  \hat \theta_0 ( \lambda).
\end{aligned}
\end{equation}

In addition, the change of variables $\lambda \rightarrow \nu,$ in the integral on $\partial D_-$ of \eqref{int_minus1} implies the transformation  $\partial D_- \rightarrow \partial D_+.$ Given that $d\nu / \Delta(\nu)= d\lambda/ \Delta (\lambda),$ we get
\begin{equation}\label{int_minus2}
- \frac1{2\pi} \int_{\partial D_-} e^{-i \lambda (L-x) - \omega (\lambda)t} \frac{\alpha_- (\lambda,t)}{\Delta (\lambda)}
  d\lambda = - \frac1{2\pi} \int_{\partial D_+} e^{-i \nu(\lambda) (L-x) - \omega (\lambda)t} \frac{\alpha_- (\nu(\lambda),t)}{\Delta (\lambda)}
  d\lambda.
\end{equation}
Thus, we can add \eqref{int_plus1} and \eqref{int_minus1} considering \eqref{int_minus2} to obtain
\begin{equation}\label{repre3}
\begin{aligned}
- \frac1{2\pi} &\int_{\partial D_+} 
 \frac{e^{- \omega (\lambda)t}}{\Delta (\lambda)} \left( e^{i\lambda x} \alpha_+ (\lambda,t) + e^{-i \nu(\lambda) (L-x)} \alpha_- (\nu(\lambda),t) \right) \\
&= - \frac1{2\pi} \int_{\partial D_+} 
 \frac{e^{- \omega (\lambda)t}}{\Delta (\lambda)} \left[ 
(2 i \lambda D_0 -K_0)  e^{-i\lambda (L-x) } +(2 i \nu(\lambda) D_0 -K_0) e^{-i \nu(\lambda) (L-x)} \right]\tilde f_0 (\omega (\lambda), t)
  d\lambda \\
&\phantom{=}+\frac1{2\pi} \int_{\partial D_+} 
 \frac{e^{- \omega (\lambda)t}}{\Delta (\lambda)} \left[ 
(2 i \lambda D_0 -K_0)  e^{i\lambda x -\tfrac{K_0}{D_0}L } \right. \\
 &\phantom{=}\left.
+(2 i \nu(\lambda) D_0 -K_0) e^{-i \nu(\lambda) (L-x)- i \lambda L}
 \right]\tilde g_0 (\omega (\lambda), t)
  d\lambda \\
  &\phantom{=}-\frac1{2\pi} \int_{\partial D_+} 
 \frac{e^{- \omega (\lambda)t}}{\Delta (\lambda)} \left[
-e^{i\lambda x} e^{-i \nu(\lambda)L} +e^{-i \nu(\lambda) (L-x)}
\right] \hat \theta_0 (\lambda)
 d \lambda \\
  &\phantom{=}-\frac1{2\pi} \int_{\partial D_+} 
 \frac{}{\Delta (\lambda)} \left[
e^{i\lambda x} e^{-i \lambda L} -e^{-i \nu(\lambda) (L-x)}
\right] \hat \theta_0 (\nu(\lambda))
 d \lambda.
\end{aligned}
\end{equation}
 We replace \eqref{repre3} into \eqref{repre2} and after some lengthy but straightforward calculations we obtain 
\begin{equation}\label{eq_final}
\begin{aligned}
\theta (x,t) &= \frac1{2\pi} \int_{\R} e^{i \lambda x - \omega (\lambda)t} \hat \theta_0 (\lambda) d\lambda \\
&\phantom{=} +\frac{i}{\pi}e^{-\frac{K_0}{2D_0}(2L-x) } \int_{\partial D_+}  \frac{e^{i \lambda L -\omega (\lambda)t} }{\Delta (\lambda)} \sin \left(x \left(\lambda+ i\frac{K_0}{2D_0} \right) \right) \hat \theta_0 (\lambda) d\lambda 
\\
&\phantom{=}+ \frac{i}{\pi}e^{-\frac{K_0}{2D_0}(L-x) } \int_{\partial D_+}  \frac{e^{- \omega (\lambda)t}}{\Delta (\lambda)} \sin \left((L-x) \left(\lambda+ i\frac{K_0}{2D_0} \right) \right) \hat \theta_0 (\nu(\lambda)) d\lambda \\
&\phantom{=}- \frac{2D_0}{\pi} e^{-\frac{K_0}{2D_0}(L-x) } \int_{\partial D_+}  \frac{e^{- \omega (\lambda)t}}{\Delta (\lambda)}\left(\lambda+ i\frac{K_0}{2D_0} \right) \sin \left((L-x) \left(\lambda+ i\frac{K_0}{2D_0} \right) \right) \tilde f_0 (\omega(\lambda),t) d\lambda \\
&\phantom{=}- \frac{2D_0}{\pi} e^{-\frac{K_0}{2D_0}(2L-x) } \int_{\partial D_+}  \frac{e^{- \omega (\lambda)t}}{\Delta (\lambda)}\left(\lambda+ i\frac{K_0}{2D_0} \right) \sin \left(x \left(\lambda+ i\frac{K_0}{2D_0} \right) \right) \tilde g_0 (\omega(\lambda),t) d\lambda .
\end{aligned}
\end{equation}

 \begin{remark}\label{rem1}
 By making the change of variables $\lambda\to\nu(\lambda)=-\lambda-i \frac{K_0}{D_0}$, we obtain a different integral representation of the solution, which reads exactly the same with \eqref{eq_final} by only substituting the contour of integration $\partial D_+$ with $\partial D_-$.
\end{remark}

\begin{remark}\label{rem2}
 Letting $L\to\infty$ and $g_0\equiv 0$, we retrieve (21) of \cite{barros19}, namely
 \begin{equation}\label{eq_final_hl}
 \begin{aligned}
\theta_{hl} (x,t) &= \frac1{2\pi} \int_{\R} e^{i \lambda x - \omega (\lambda)t} \hat \theta_0 (\lambda) d\lambda \\
&\phantom{=}- \frac{1}{2\pi} \int_{\partial D_+}  
e^{i \lambda x - \omega (\lambda)t} \left[ \hat \theta_0 (\nu(\lambda)) + \left(2i\lambda D_0- K_0 \right)  \tilde f_0 (\omega(\lambda),t) \right] d\lambda.
\end{aligned}
 \end{equation}
\end{remark}

\section{Test examples}\label{sec_examples}

We examine the feasibility of the derived formula in three different test cases modelling the following infiltration processes: flooding (Example 1), rainfall with constant flux (Example 2) and rainfall at dry soil with a water tank at the bottom surface (Example 3). 

\subsection*{Example 1: Flooding}

In the first case, the
 initial and boundary conditions of \eqref{bvp} are simplified as  
\begin{subequations}\label{bc1}
\begin{alignat}{3}
\theta (x,0) &= \theta_0, \quad && 0 <x <L, \\ 
\theta (0,t) &= \theta_1, \quad &&t>0, \\
 \theta (L,t) &= \theta_0, \quad &&t>0, 
\end{alignat}
\end{subequations}
for $\theta_1 > \theta_0 >0,$ constant values. Without loss of generality, we consider an equivalent problem for $u(x,t) = \theta(x,t) - \theta_0,$ satisfying 
\begin{subequations}\label{bc1u}
\begin{alignat}{3}
\frac{\partial u}{\partial t} + K_0 \frac{\partial u}{\partial x} &= D_0 \frac{\partial^2 u}{\partial x^2},  \quad && 0<x <L, \, t>0,  \\
u (x,0) &= 0, \quad && 0 <x <L, \\ 
u (0,t) &= \theta_1 - \theta_0, \quad &&t>0, \label{bc1u3}\\
u (L,t) &= 0, \quad &&t>0. 
\end{alignat}
\end{subequations}

Given \eqref{bc1u3}, we get
\begin{equation}
\tilde f_0 (\omega(\lambda),t) = \int_0^t e^{\omega(\lambda) \tau} f(\tau )d\tau
=  (\theta_1 - \theta_0)  \int_0^t e^{\omega(\lambda) \tau} d\tau =  (\theta_1 - \theta_0) \frac{e^{\omega(\lambda) t}-1}{\omega(\lambda)}
\end{equation}
and using that $\hat \theta_0 (\lambda) = \hat \theta_0 (\nu(\lambda)) = \tilde g_0 (\omega(\lambda)) \equiv 0,$ the solution of \eqref{bc1u} considering \eqref{eq_final} is given by
\begin{equation}\label{sol_ex1}
\begin{aligned}
u_1 (x,t) = (\theta_0 - \theta_1) \frac{2D_0}{\pi} e^{-\tfrac{K_0}{2D_0} (L-x)} \int_{\partial D_+} \frac{1 - e^{-\omega(\lambda)t}}{\omega (\lambda) \Delta (\lambda)} \left(\lambda + i \frac{K_0}{2D_0} \right) \sin \left( (L-x) \left(\lambda + i \frac{K_0}{2D_0} \right)\right) d\lambda .
\end{aligned}
\end{equation}
The above expression can readily be evaluated numerically; however, in what follows we provide some analytical considerations which yield a more efficient form of the solution from computational aspect. 

Solving $\omega (\lambda) \Delta (\lambda) =0,$ we get the two isolated roots $\lambda=0$ and $\lambda=-i \tfrac{K_0}{D_0}$, as well as the family of roots $\lambda_n=-i \tfrac{K_0}{2D_0} + \frac{n \pi}{L}, \ n\in\Z$. It is obvious that $\lambda=0$ is a removable singularity of the integrand, and the other roots reside on the lower half complex $\lambda$-plane. Thus, we are allowed to compute the integral of \eqref{sol_ex1} along the slightly deformed contour $\partial \tilde{D}_+$, which is defined as $\partial D_+$, but it bypasses $\lambda=0$ with a circular arc of radius $\epsilon\to0$.

Now we are able to compute the integral of \eqref{sol_ex1} by decomposing it as $I= I_1 + I_2,$ where
\begin{equation}\label{int_decomp}
\begin{aligned}
I_1 (x,t) &= \int_{\partial \tilde{D}_+} \frac{1}{\omega (\lambda) \Delta (\lambda)} \left(\lambda + i \frac{K_0}{2D_0} \right) \sin \left( (L-x) \left(\lambda + i \frac{K_0}{2D_0} \right)\right) d\lambda, \\
I_2 (x,t) &=  -\int_{\partial \tilde{D}_+} \frac{ e^{-\omega(\lambda)t}}{\omega (\lambda) \Delta (\lambda)} \left(\lambda + i \frac{K_0}{2D_0} \right) \sin \left( (L-x) \left(\lambda + i \frac{K_0}{2D_0} \right)\right) d\lambda.
\end{aligned}.
\end{equation}
Since the integrand of $I_1$ is bounded and analytic on $D_+$, then Cauchy's residue theorem implies that the contribution of $I_1$ is given as $1/4$ of the residue contribution at  $\lambda=0$, namely
\[
I_1 (x,t) =- \frac{\pi}{D}\frac{\sinh \left(\frac{K_0}{2D_0} (L-x) \right)}{1-e^{-\tfrac{K_0 L}{D_0}}}.
\]
The second integral using change of variables can be also further simplified but it can also be directly numerically computed: this simplification does not provide a significant computational improvement.

\subsection*{Example 2: Rainfall - constant flux}

In this case, the initial and boundary conditions are written as
\begin{subequations}\label{bc2}
\begin{alignat}{3}
\theta (x,0) &= \theta_0, \quad && 0 <x <L, \\ 
\theta (0,t) &= f_2(t), \quad &&t>0, \\
 \theta (L,t) &= \theta_0, \quad &&t>0. 
\end{alignat}
\end{subequations}

As in the first example, we define $u(x,t) = \theta(x,t) - \theta_0,$ and we consider instead the problem 
\begin{subequations}\label{bc2u}
\begin{alignat}{3}
\frac{\partial u}{\partial t} + K_0 \frac{\partial u}{\partial x} &= D_0 \frac{\partial^2 u}{\partial x^2},  \quad && 0<x <L, \, t>0,  \\
u (x,0) &= 0, \quad && 0 <x <L, \\ 
u (0,t) &= f_2(t) - \theta_0, \quad &&t>0, \label{bc2u2}\\
u (L,t) &= 0, \quad &&t>0.
\end{alignat}
\end{subequations}

The boundary equation \eqref{bc2u2} now gives
\begin{equation}
\tilde f_0 (\omega(\lambda),t) = \int_0^t e^{\omega(\lambda) \tau} (f_2(\tau ) - \theta_0 )d\tau
= \tilde f_2 (\omega(\lambda),t)
-\theta_0 \frac{e^{\omega(\lambda) t}-1}{\omega(\lambda)}.
\end{equation}
Thus, the solution $u_2$ of this problem, considering \eqref{sol_ex1}, takes the form
\begin{equation}
\begin{aligned}
u_2(x,t) &= \frac{\theta_0}{\theta_0 - \theta_1} u_1 (x,t) \\
&\phantom{=}- 
 \frac{2D_0}{\pi} e^{-\tfrac{K_0}{2D_0} (L-x)} \int_{\partial D_+} \frac{ e^{-\omega(\lambda)t}}{ \Delta (\lambda)} \tilde f_2 (\omega(\lambda),t) \left(\lambda + i \frac{K_0}{2D_0} \right) \sin \left( (L-x) \left(\lambda + i \frac{K_0}{2D_0} \right)\right) d\lambda.
\end{aligned}
\end{equation}

The numerical evaluation of the last term in the above equation depends on the specific form of $f_2,$ which will be specified in the next section. 

\subsection*{Example 3: Rainfall - water tank}

Starting from \eqref{eq_rich} one can derive a PDE for the water pressure head $h$. We define $h_{a,\epsilon} = e^{a h}-\epsilon.$ Then, following \cite[Eqs. (19) -- (21)]{Tracy11}, the IBVP for the transformed pressure head reads
\begin{subequations}\label{ex3}
\begin{alignat}{3}
\frac{\partial h_{a,\epsilon}}{\partial t} - K_1 \frac{\partial h_{a,\epsilon}}{\partial x} &= D_1 \frac{\partial^2 h_{a,\epsilon}}{\partial x^2},  \quad && 0<x <L, \, t>0,  \\
h_{a,\epsilon} (x,0) &= 0, \quad && 0 <x <L, \\ 
h_{a,\epsilon} (0,t) &= 0, \quad &&t>0, \\
h_{a,\epsilon} (L,t) &= 1-\epsilon, \quad &&t>0,  \label{bc33}
\end{alignat}
\end{subequations}
for some positive constant coefficients $K_1, \, D_1$ and $\epsilon.$ The main difference compared to the previous two examples is the appearance of the negative sign in front of the first spatial partial derivative. This is because the unknown function is now the pressure head and not the water content.   We apply the same formulas from before assuming $K_1 = - K_0 <0.$

Since now the only non-zero boundary condition is the one at $x=L,$ we consider \eqref{bc33} and we get
\begin{equation}
\tilde g_0 (\omega(\lambda),t) = \int_0^t e^{\omega(\lambda) \tau} g(\tau )d\tau
=  (1-\epsilon) \frac{e^{\omega(\lambda) t}-1}{\omega(\lambda)}
\end{equation}
and the solution, using \eqref{eq_final} as well as Remark \ref{rem1}, is given by
\begin{equation}\label{sol_ex3}
\begin{aligned}
h_{a,\epsilon} (x,t) = (\epsilon-1) \frac{2D_1}{\pi} e^{-\tfrac{K_1}{2D_1} (2L-x)} \int_{\partial D_-} \frac{1 - e^{-\omega(\lambda)t}}{\omega (\lambda) \Delta (\lambda)} \left(\lambda + i \frac{K_1}{2D_1} \right) \sin \left( x\left(\lambda + i \frac{K_1}{2D_1} \right)\right) d\lambda
\end{aligned}
\end{equation}
The computation of this integral follows the exact same steps as the one in Example 1.



\section{Numerical Implementation}\label{sec_numerics}

We compute numerically the integrals over $\partial D_+$ by considering the parametrized half-line $\gamma (\theta,t) = e^{i\theta} t,$ for $t\geq 0,$ where the angle $\theta$ is chosen such that the contour $\gamma$ lies between the real axis and $\partial D_+.$ Thus, we get
\begin{align*}
\int_{\partial D_+} \phi (\lambda) d\lambda &= e^{i\theta} \int_0^\infty \phi (\gamma (\theta,t)) dt -  e^{i (\pi-\theta)} \int_0^\infty \phi (\gamma (\pi -\theta,t)) dt\\
&= \int_0^\infty e^{i\theta}\phi \left(e^{i\theta} t\right) + e^{-i \theta}\phi \left(-e^{-i \theta} t\right) dt
\end{align*}

\subsection*{Example 1: Flooding}

\subsubsection*{Case 1}
\begin{figure}
    \centering
    \begin{subfigure}[t]{0.5\textwidth}
        \centering
        \includegraphics[width=0.9\textwidth]{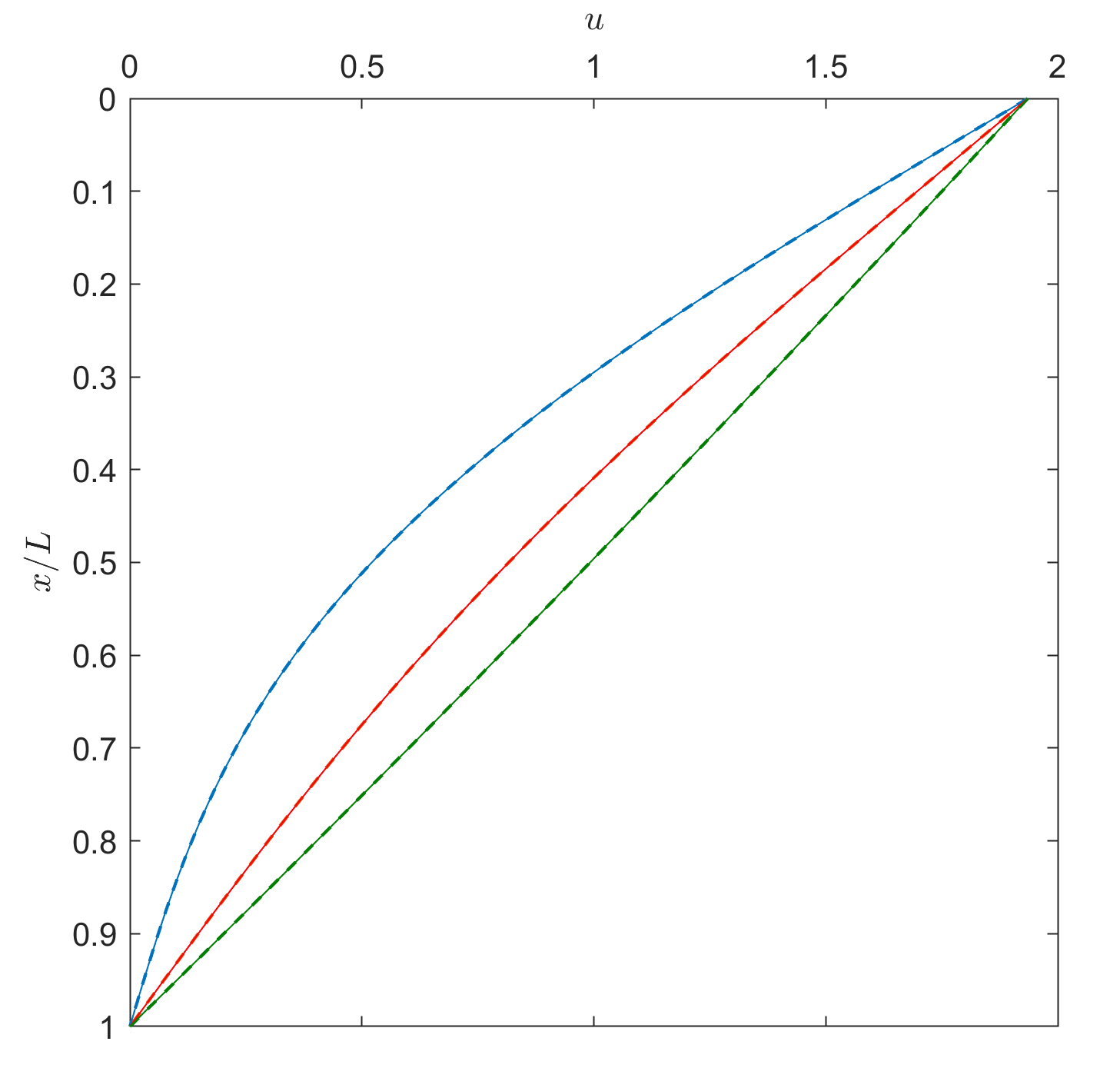}
    \end{subfigure}%
    ~ 
    \begin{subfigure}[t]{0.5\textwidth}
        \centering
        \includegraphics[width=0.9\textwidth]{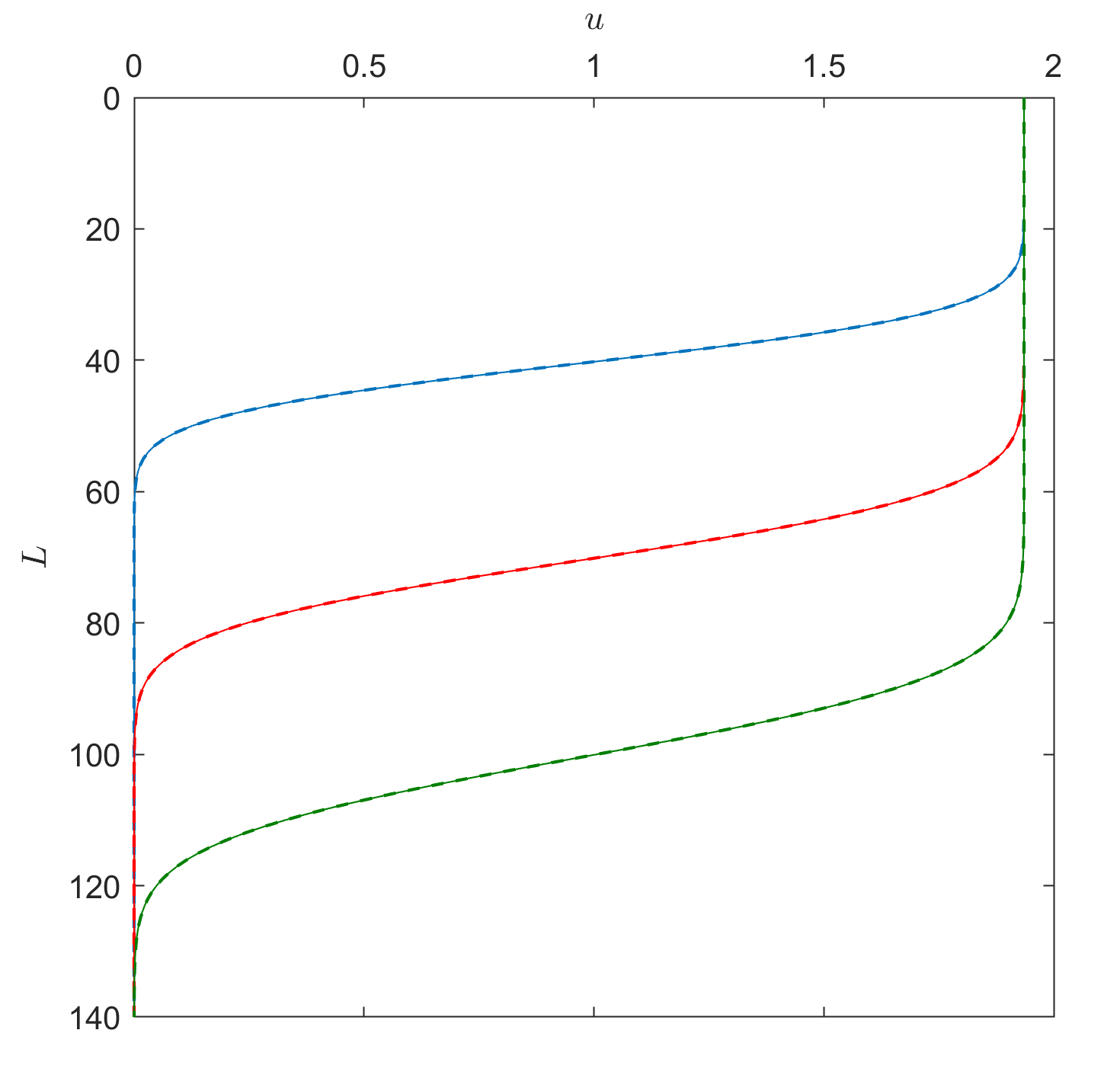}
    \end{subfigure}
    \caption{Left: The solutions $u_1$ (solid line) and $u_P$ (dotted line) for $L = 0.05 \si{\cm}$ for the first case of the first example, at the time steps: $t_1 = 0.03 \si{\sec}$ (blue line), $t_2 = 0.06 \si{\sec}$ (red line) and  $t_3 = 0.6 \si{\sec}$ (green line).    Right:  The solutions $u_1$ (solid line) and $u_P$ (dotted line) for $L = 140 \si{\cm},$ at $t_1 = 40 \si{\min}$ (blue line), $t_2 = 70 \si{\min}$ (red line) and  $t_3 = 100 \si{\min}$ (green line). }\label{fig1a}
\end{figure}

In the special case of $D_0 = 1/2$ and $K_0 = 1,$ Philip in \cite{Philip91} presented a solution suitable for shallow profiles. We compare our solution $u_1$ with the series representation \cite[Eq. (27)]{Philip91}, to be called $u_P.$ We set $L=0.05$ and $\theta_1 - \theta_0 = 1.9355$ in order to reproduce  \cite[Fig. 3]{Philip91}. In the left picture of \autoref{fig1a} we observe that the solutions coincide for the different time steps.

As our solution remains valid for all $L>0$, we keep the same setup as before and investigate the convergence of the two solutions over extended profiles and longer durations. It is noteworthy that achieving the level of convergence depicted in the right in the right picture of \autoref{fig1a}, we had to consider a greater number of terms in the series representation (from tens to hundreds), resulting in a computationally expensive process.

\subsubsection*{Case 2}

 We set $\theta_0 = 0.025 \si[per-mode=symbol]
{\cm^3\per\cm^3}$ and  $\theta_1 = 0.335 \si[per-mode=symbol]
{\cm^3\per\cm^3}$ as initial and surface soil water content, respectively.  The conductivity is given by $K_0 = 4.32 \si[per-mode=symbol]
{\cm\per\hour}$ and the diffusivity by $D_0 = 0.4653 \si[per-mode=symbol]
{\cm^2\per\second}$ \cite{Arg97}.

We consider different soil thicknesses  to show that the larger the thickness the better the approximation of the widely used solution derived by Philip \cite{philip66}
\begin{equation}
\theta_P (x,t) = \frac{1}{2} \left[\mbox{Erfc}\left(\frac{x-\kappa t}{2\sqrt{D_0t}}\right) + e^{\tfrac{\kappa x}{D_0}} \mbox{Erfc}\left(\frac{x+\kappa t}{2\sqrt{D_0t}}\right)\right],  \quad \kappa = \frac{K_0}{\theta_1 - \theta_0}. 
\end{equation}

\begin{figure}
    \centering
    \begin{subfigure}[t]{0.5\textwidth}
        \centering
        \includegraphics[width=0.9\textwidth]{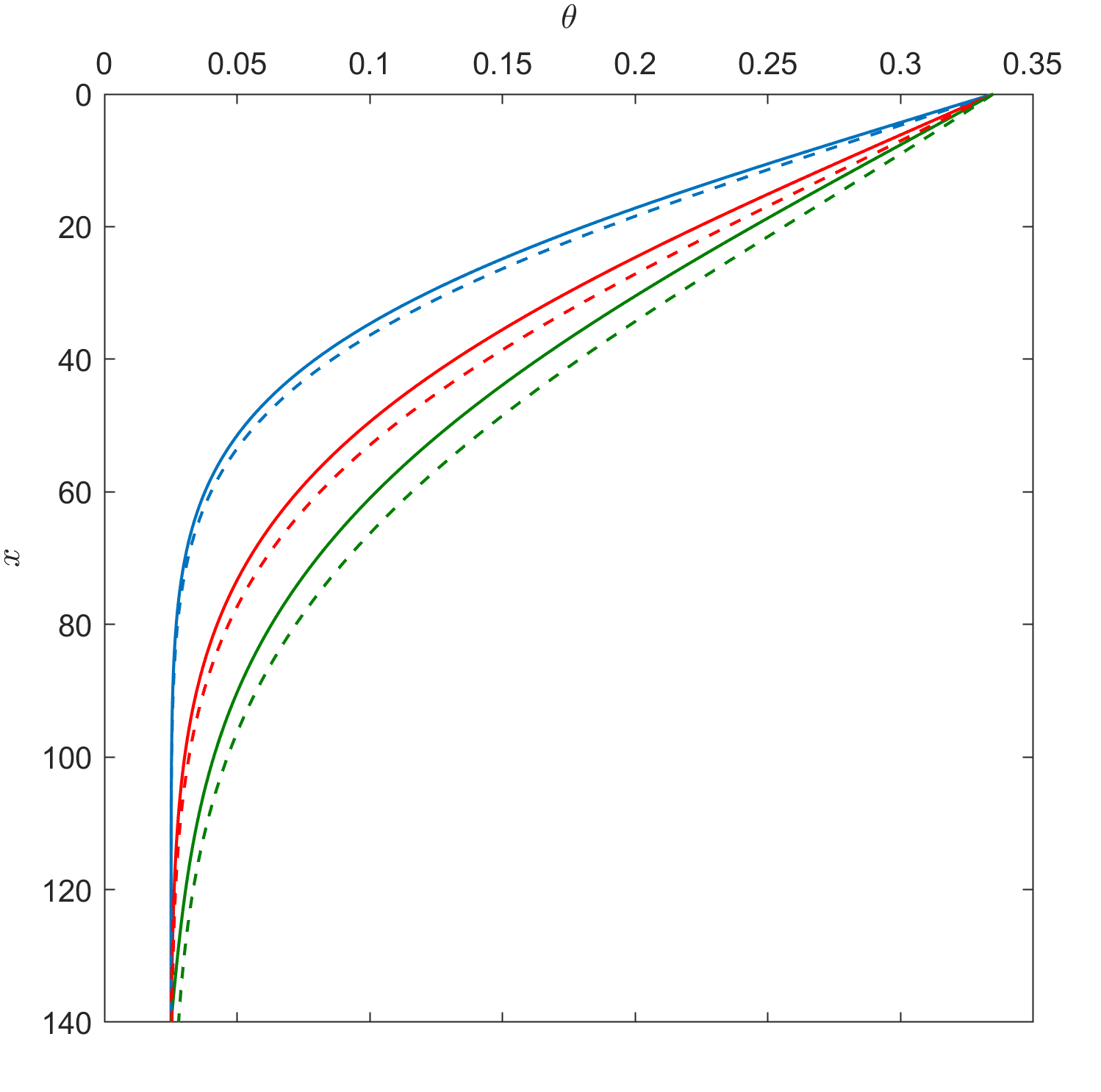}
    \end{subfigure}%
    ~ 
    \begin{subfigure}[t]{0.5\textwidth}
        \centering
        \includegraphics[width=0.9\textwidth]{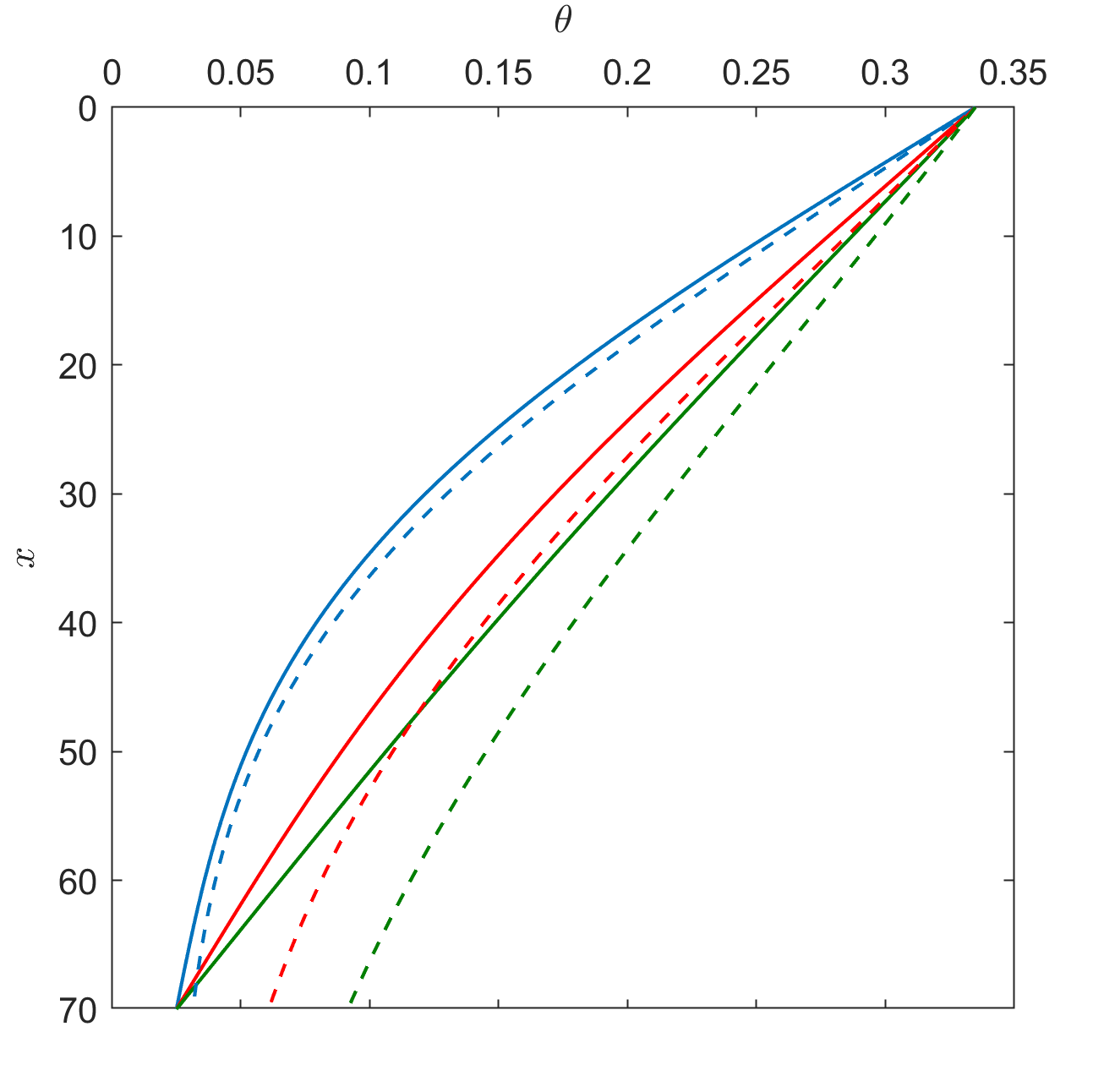}
    \end{subfigure}
    \caption{The solutions $\theta$ (solid line) and $\theta_P$ (dotted line) for $L = 140 \si{\cm}$ (left picture) and $L = 70 \si{\cm}$ (right picture) of the first example (case 2). The different colors indicate the different time steps: $t_1 = 15 \si{\min}$ (blue line), $t_2 = 30 \si{\min}$ (red line) and  $t_3 = 45 \si{\min}$ (green line). }\label{fig1}
\end{figure}

In \autoref{fig1} we present the two solutions $\theta$ (solid line) and $\theta_P$ (dotted line) for $L = 140 \si{\cm}$ (left picture) and $L = 70 \si{\cm}$ (right picture). The two solutions are compared at the times $t_1 = 15 \si{\min}$ (blue line), $t_2 = 30 \si{\min}$ (red line) and at $t_3 = 45 \si{\min}$ (green line). 
As discussed in \autoref{intro}, we observe that the solution of \eqref{bvp} can be approximated well by $\theta_P$ for large $L$ and/or small times $t.$ However, if $L$ is not large enough the solution displays different features which $\theta_P$ fails to capture.

The solution $\theta_P$ was also compared with the one derived using Fokas method for the half-line \cite[Eq. (25)]{barros19}. We get a perfect match with that one for large soil thickness.


\subsection*{Example 2: Rainfall - constant flux}

For this example we compare our solution with the one derived by Braester in \cite{Bra73} for the special case of $K_0 = a D_0,$ for some positive constant $a.$
 We keep all parameters the same as in the first example, except for the diffusivity $D_0 = 0.387 \si[per-mode=symbol]
{\cm^2\per\second}.$ 

The boundary function, defined by $f_2,$ is given 
\[
f_2 (t) = \frac12  \left[\mbox{Erfc}\left(\frac{- \sqrt{K_a t}}{2}\right) - (1+ K_a t) \mbox{Erfc}\left(\frac{ \sqrt{K_a t}}{2}\right) + 2 \sqrt{\frac{K_a t}{\pi}} e^{-\tfrac{K_a t}4{}}\right], \quad K_a = a K_0. 
\]
We define by $\theta_B$ the corresponding solution for the half line, see \cite[Eq. (12)]{Bra73}. As in the previous example, we choose a relative large $L$ so that the comparison makes sense.

\begin{figure}
    \centering
    \begin{subfigure}[t]{0.5\textwidth}
        \centering
        \includegraphics[width=0.89\textwidth]{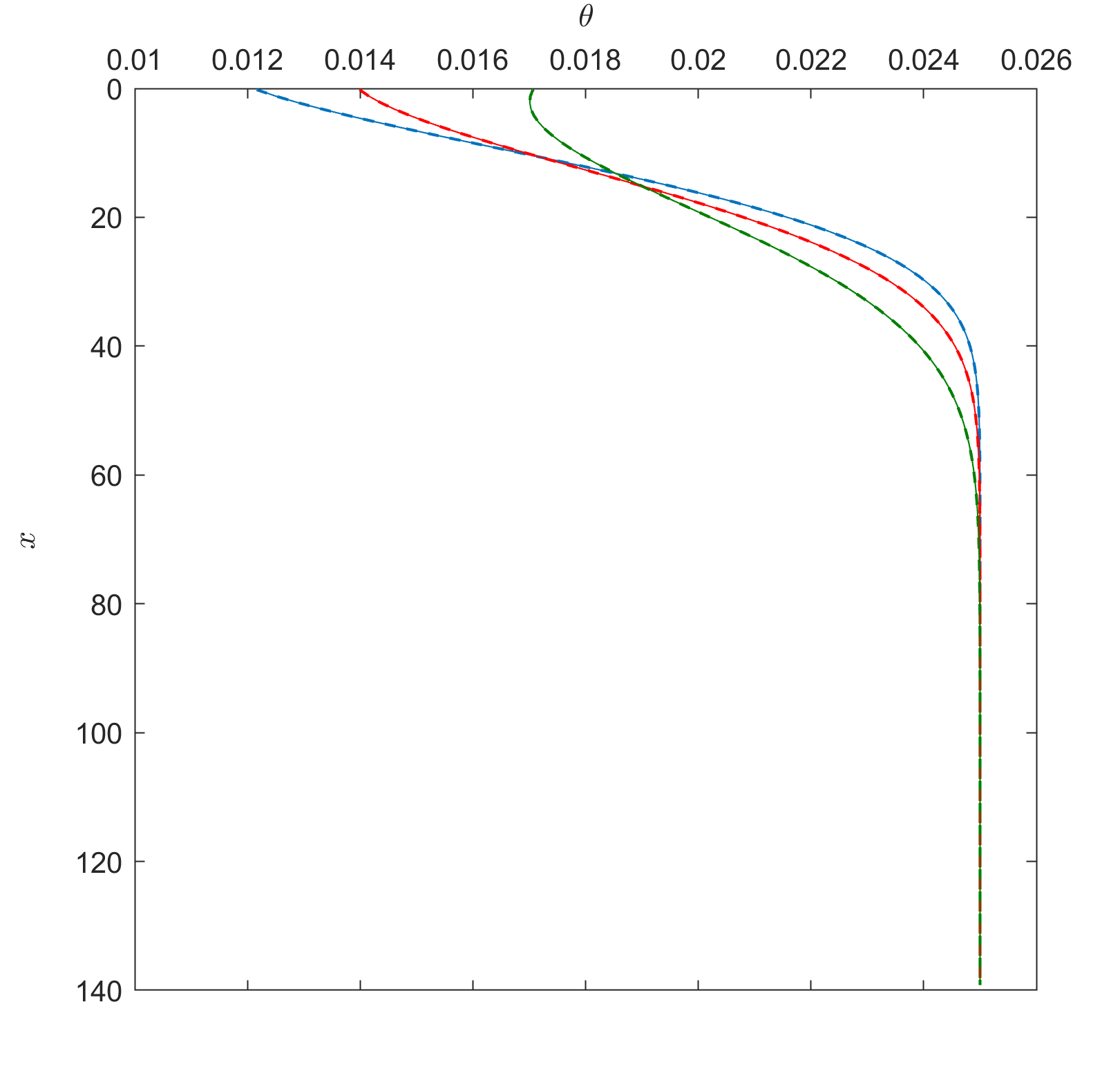}
    \end{subfigure}%
    ~ 
    \begin{subfigure}[t]{0.5\textwidth}
        \centering
       \includegraphics[width=0.9\textwidth]{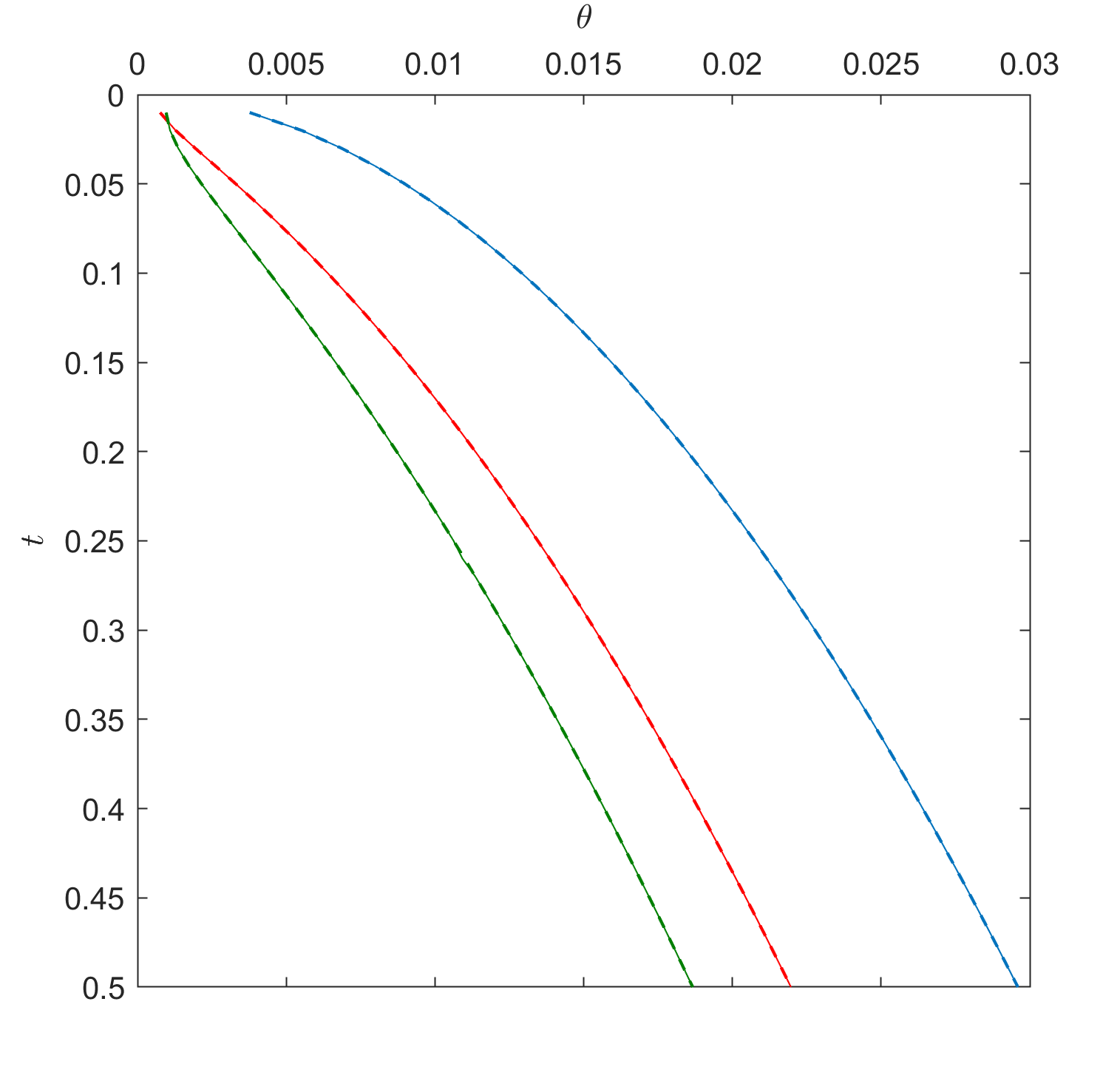}
    \end{subfigure}
    \caption{The solutions $\theta$ (solid line) and $\theta_B$ (dotted line) for $L = 140 \si{\m},$ of the second example. Left: The solutions at $t_1 = 5 \si{min}$ (blue line) $t_2 = 7 \si{min}$ (red line) and $t_3 = 10 \si{min}$ (green line). Right: The solutions  at $x_1= 0.5 \si{m}$ (blue), $x_2 =10\si{m}$ (red) and $x_3 =15 \si{m}$ (green) for time in hours.  
     }\label{fig2}
\end{figure}

The plots of the functions are presented in \autoref{fig2} for  $L = 140 \si{\m}$ and $a = 0.001.$
The solutions are compared at the time steps: $t_1 = 5 \si{min}$ (blue line) $t_2 = 7 \si{min}$ (red line) and $t_3 = 10 \si{min}$ (green line) in the left picture and at the positions $x_1= 0.5 \si{m}$ (blue), $x_2 =10\si{m}$ (red) and $x_3 =15 \si{m}$ (green) for $t\in [0\si{min},\,30\si{min}],$ in the right picture.

\subsection*{Example 3: Rainfall - water tank}

We set 
\[
K_1 = -\frac{K_s}{\theta_1 - \theta_0}, \quad D_1 = \frac{K_s}{a(\theta_1 - \theta_0)},
\]
where $K_s$ is the hydraulic conductivity of the saturated soil.  The solution derived by Tracy \cite{Tracy11} admits the form
\begin{equation}\label{tracy1}
h_T (x,t) = \frac{1}{a} \log\left[ (1-\epsilon) e^{\tfrac{a}{2} (L-x)} \left( \frac{\mbox{sinh} \left(\frac{a x}{2} \right)}{\mbox{sinh} \left(\frac{a L}{2} \right)} + \frac{2 D_1}{L} \sum_{k=1}^\infty (-1)^k \frac{\lambda_k}{\mu_k} \sin (\lambda_k x) e^{- \mu_k t} \right) +\epsilon \right],
\end{equation}
where $\lambda_k =  \tfrac{\pi k}{L}$ and $\mu_k = D_1 (\tfrac{a^2}{4} + \lambda_k^2).$ In the following two cases, we compare $h$ as a function of $h_{a,\epsilon},$ given by \eqref{sol_ex3} with $h_T$ given by \eqref{tracy1}. 

\subsubsection*{Case 1}

We consider the parameters used in \cite{Tracy11}, meaning the pressure head for dry soil is $h_0 = -20 \si{m},$ $a = 0.1 \si{m}^{-1},$ $L=50 \si{m},$ and $K_s = 0.1 \times 10^{-5} \si[per-mode=symbol]{m\per\day}.$ The water contents are given by $\theta_1 = 0.45$ and $\theta_0 = 0.15.$

\begin{figure}
    \centering
    \begin{subfigure}[t]{0.5\textwidth}
        \centering
        \includegraphics[width=0.87\textwidth]{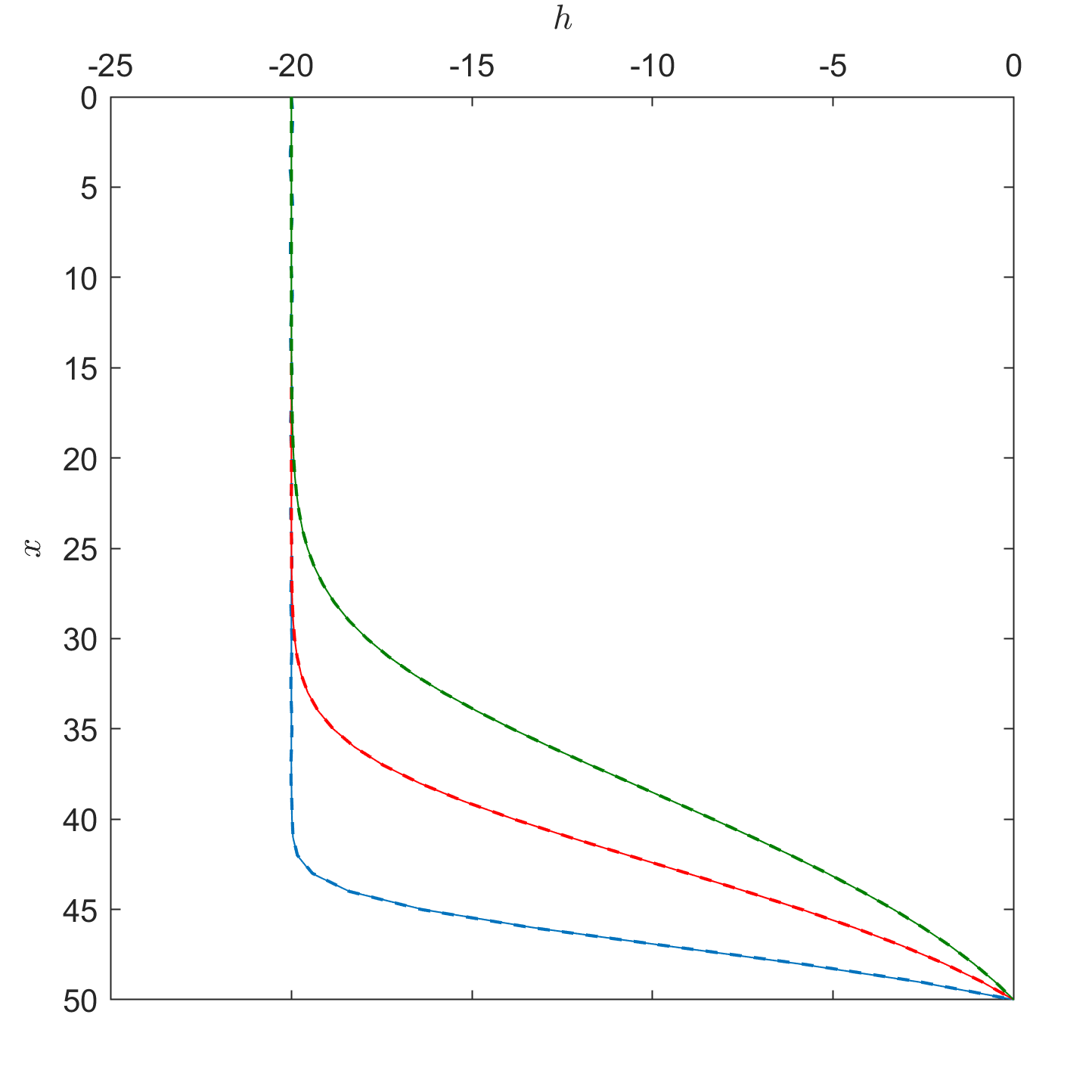}
    \end{subfigure}%
    ~ 
    \begin{subfigure}[t]{0.5\textwidth}
        \centering
        \includegraphics[width=0.9\textwidth]{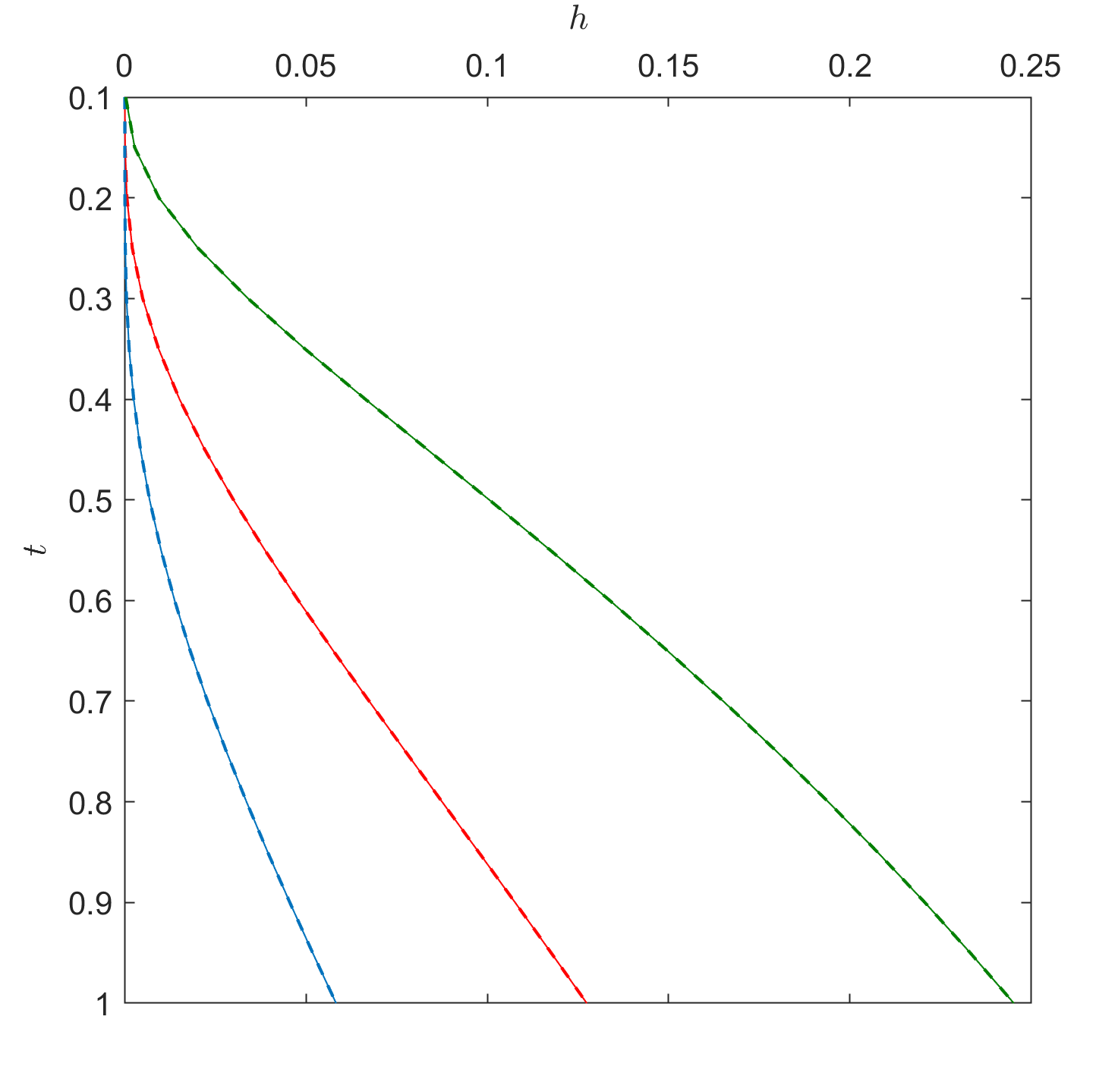}
    \end{subfigure}
    \caption{Case 1 of 3rd example: The solution $h$ (solid line) and $h_T$ (dotted line) for $L = 50 \si{\m}.$ Left: The solutions at $t_1 = 1\si{day}$ (blue line) $t_2=5 \si{day}$ (red line) and $t_3=10 \si{day}$ (green line). Right: The solutions  at $x_1=45\si{m}$ (blue), $x_2=46\si{m}$ (red) and $x_3=47\si{m}$ (green).  
     }\label{fig4}
\end{figure}

The two solutions are compared once at different time steps (left picture) and then at different positions (right) in \autoref{fig4}. In the left picture, our solution (solid line) is plotted against $h_T$ (dotted line) at the time steps $t_1 = 1\si{day},$ (blue line) $t_2= 5\si{day}$ (red line)  and at $t_3 = 10\si{day}$ (green line). In the right picture, the two solutions are presented at the positions $x=45\si{m}$ (blue),  $x=46\si{m}$ (red) and $x=47 \si{m}$ (green) for $t\in [0.1\si{day},\,1\si{day}]. $

\subsubsection*{Case 2}

\begin{figure}
    \centering
    \begin{subfigure}[t]{0.5\textwidth}
        \centering
        \includegraphics[width=0.9\textwidth]{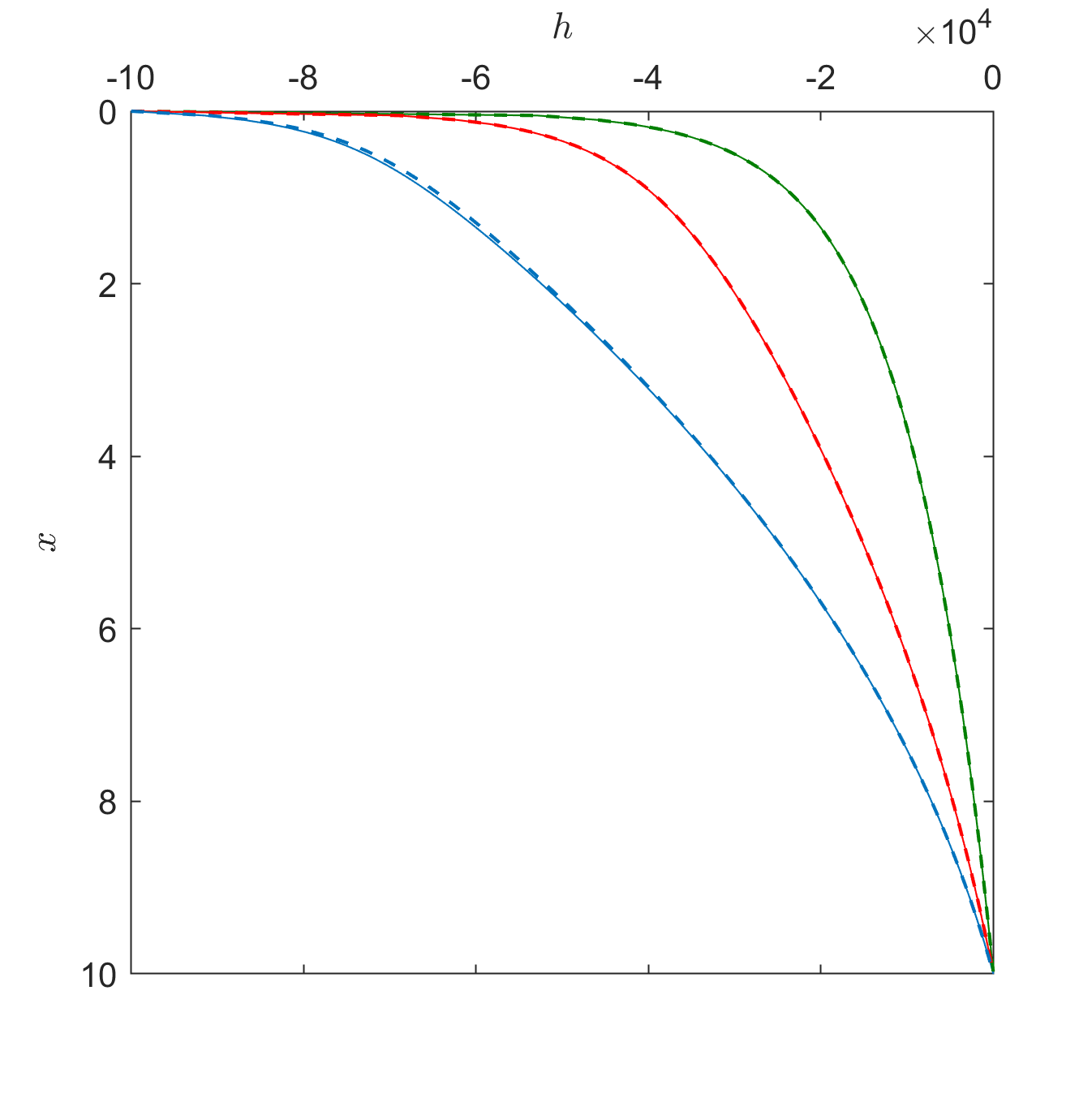}
    \end{subfigure}%
    ~ 
    \begin{subfigure}[t]{0.5\textwidth}
        \centering
        \includegraphics[width=0.9\textwidth]{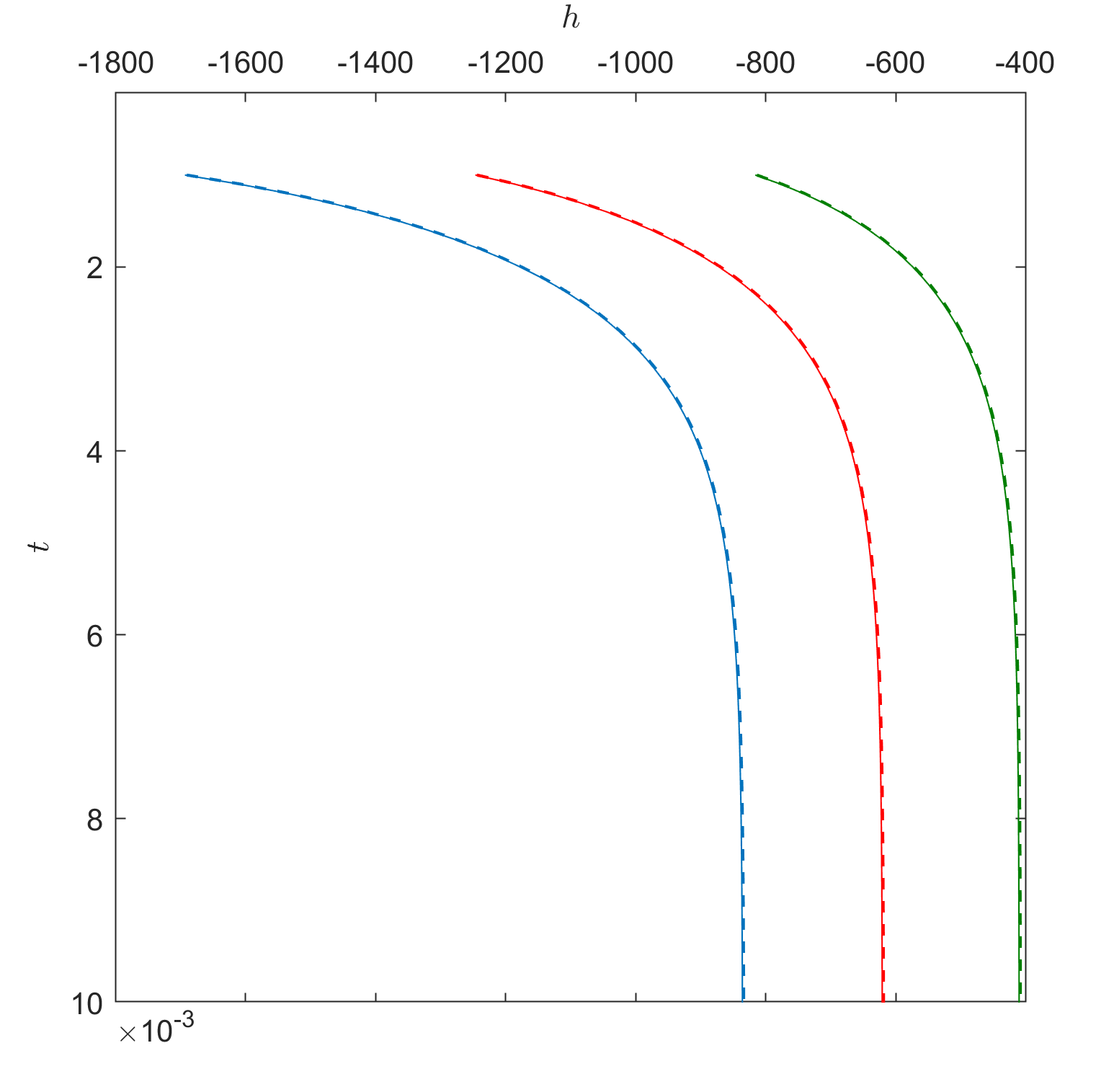}
    \end{subfigure}
    \caption{Case 2 of 3rd example: The solution $h$ (solid line) and $h_T$ (dotted line) for $L = 10 \si{\m}.$ Left: The solutions at $t_1 = 0.0005 \si{h}$ (blue line) $t_2=0.001 \si{h}$ (red line) and $t_3=0.01 \si{h}$ (green line). Right: The solutions  at $x_1=9.2\si{m}$ (blue), $x=9.4\si{m}$ (red) and $x=9.6\si{m}$ (green).  
     }\label{fig3}
\end{figure}

 We set $L = 10 \si{m}$ and $h_0 = -10^5,$ the pressure head when the soil is dry, $a= 10^{-4},$ $K_s = 9 \times 10^{-5} \si[per-mode=symbol]
{m\per\second},$ and the water contents are given by $\theta_1 = 0.5$ and $\theta_0 = 0.11$ \cite{liu17}. In the left picture of \autoref{fig3}, we see our solution (solid line) compared with $h_T$ (dotted line) at $t_1 = 0.0005\si{h},$ (blue line) $t_2= 0.001\si{h}$ (red line)  and at $t_3 = 0.01 \si{h}$ (green line). In the right picture, we present the matching between the two solutions at the specific positions $x=9.2\si{m}$ (blue),  $x=9.4\si{m}$ (red) and $x=9.6\si{m}$ (green) for $t\in [10^{-3}\si{h},\,10^{-2}\si{h}]. $

\section*{Conclusions}

In this work we reviewed and enriched the mathematical modelling of the vertical infiltration in bounded profiles, by formulating the initial boundary value problem \eqref{bvp}, to be considered in the finite interval $(0,L)$; for example, the boundary condition which fixes the water content at the bottom of the interval $x=L$ can be seen as the physical set-up which does not allow for any infiltration at this point.  We derived the solution of the latter problem in the form of the integral representation \eqref{eq_final}, which provides the solution for general initial and boundary conditions in a rather effective and computationally cheap form due to the exponential decay of the relevant integrands.  We note that the limit $L\to\infty$ yields \eqref{eq_final_hl}, which is  the solution for the half-line problem, and coincides with \cite[Equation (21)]{barros19}.
This analytical result reflects the physical fact that if the length $L$ of the interval is large enough then the classical solutions, e.g. \cite{philip1957}, yield sufficient results for the problem \eqref{bvp}. However, if the length $L$ is small, then the solution \eqref{eq_final} displays different features both quantitatively and qualitatively, which cannot be captured by modifying  solutions given for the half-line problem. 
 
This work can be seen as a first step to apply the Fokas method to problems which model infiltration under various conditions. We compare our analytical results, which have been derived in a unified way, to classical results which have been extracted by using various methodologies. This is known as the direct problem and as future work we aim to consider the corresponding inverse problem in two directions: first, given boundary measurements estimate the soil parameters, namely diffusivity and conductivity; second, control the water distribution throughout the finite interval by applying appropriate boundary conditions at the bottom $x=L$. Considering the latter inverse problem, the methodology introduced by one of the authors \cite{KOD23}, suggests an effective computational algorithm; already this general methodology was applied in \cite{hwang23} with negligible modifications, for controlling a specific case of the \eqref{bvp}. 

\subsection*{Declarations}

\subsubsection*{Author Contributions} Formal analysis: I.A. and L.M.; Investigation: K.K. and L.M.; Methodology: K.K. and L.M.; Numerical results: I.A. and L.M.; 
Writing-original draft: all authors. 

\subsubsection*{Data availability statement} The manuscript has no associated data.

\subsubsection*{Conflict of interest} The authors declare no conflict of interest.

\bibliographystyle{siam}
\bibliography{refs}

\end{document}